\documentclass[11pt, a4paper, reqno]{amsart}
\usepackage{preambolo}

\begin{document}


\setcounter{secnumdepth}{3}

\setcounter{tocdepth}{2}

\title{\textbf{On the Kodaira-Spencer's Problem on almost Hermitian $4$-manifolds}}

\author[Lorenzo Sillari]{Lorenzo Sillari}

\address{Lorenzo Sillari: Scuola Internazionale Superiore di Studi Avanzati (SISSA), Via Bonomea 265, 34136 Trieste, Italy.} 
\email{lsillari@sissa.it}

\author[Adriano Tomassini]{Adriano Tomassini}

\address{Adriano Tomassini: Dipartimento di Scienze Matematiche, Fisiche e Informatiche, Unità di Matematica e
Informatica, Università degli Studi di Parma, Parco Area delle Scienze 53/A, 43124, Parma, Italy}
\email{adriano.tomassini@unipr.it}

\maketitle

\begin{abstract} 
\noindent \textsc{Abstract}. 
In 1954, Hirzebruch reported a problem posed by Kodaira and Spencer: on compact almost complex manifolds, is the dimension $h^{p,q}_{\bar \partial}$ of the kernel of the Dolbeault Laplacian independent of the choice of almost Hermitian metric? In this paper, we review recent progresses on the original problem and we introduce a similar one: on compact almost complex manifolds, find a generalization of Bott-Chern and Aeppli numbers which is metric-independent. We find a solution to our problem valid on almost K\"ahler $4$-manifolds.

\end{abstract}

\blfootnote{  \hspace{-0.55cm} 
{\scriptsize 2020 \textit{Mathematics Subject Classification}. Primary: 32Q60, 58A14; Secondary: 53C15. 

\noindent The authors are partially supported by GNSAGA of INdAM. The second author is partially supported by the Project PRIN 2017 “Real and Complex Manifolds: Topology, Geometry and
holomorphic dynamics” (code 2017JZ2SW5).}}

\section{Introduction}\label{sec:intro}

The aim of this paper is to formulate a problem similar to that posed by Kodaira and Spencer and appeared in Hirzebruch's 1954 list of problems \cite[Problem 20]{Hir54}, and to find a solution to our problem valid on almost K\"ahler $4$-manifolds. Kodaira and Spencer asked to generalize \emph{Hodge numbers}, usually defined for complex manifolds, to almost complex manifolds. Our problem, as we will see later, consists in finding a generalization of \emph{Bott-Chern and Aeppli numbers}.\\
More precisely, let $(M,J)$ be a compact almost complex manifold. Fix an \emph{almost Hermitian metric}, that is a $J$-compatible metric inducing a Hermitian metric on the fibers of $TM$. As observed by Hirzebruch \cite{Hir54}, the \emph{Dolbeault Laplacian} 
\[
\Delta_{\bar \partial} := \bar \partial \bar \partial  ^* + \bar \partial  ^* \bar \partial  
\]
is a self-adjoint elliptic operator, independently of whether $J$ is a complex structure or not, and its kernel computed on $(p,q)$-forms, the so called \emph{space of $\bar \partial$-harmonic $(p,q)$-forms}, is a finite-dimensional vector space over $\C$ of complex dimension $h^{p,q}_{\bar \partial}$. In the complex case, $h^{p,q}_{\bar \partial}$ is independent of the choice of Hermitian metric by Hodge theory.\\
The following question appeared as Problem 20 in Hirzebruch's 1954 list of problems.
\vspace{.2cm}

\textbf{Kodaira-Spencer's Problem.} (Kodaira-Spencer, see \cite[Problem 20]{Hir54})\\
Is $h^{p,q}_{\bar \partial}$ independent of the choice of (almost) Hermitian structure? If not, give some other definition of the $h^{p,q}_{\bar \partial}$ which depends only on the almost complex structure and which generalizes the $h^{p,q}_{\bar \partial}$ of a complex manifold.
\vspace{.2cm}

Recently, Holt and Zhang \cite{HZ22a} answered negatively to Kodaira-Spencer's question, and several authors gave various generalizations of $h^{p,q}_{\bar \partial}$ to almost complex manifolds (see section \ref{subsec:prel:KS} for a detailed review).\\
Similarly to what happens for the Dolbeault Laplacian, the \emph{Bott-Chern Laplacian}
\[
\Delta_{\partial + \bar \partial} := \partial \bar \partial ( \partial \bar \partial)^* + ( \partial \bar \partial)^* \partial \bar \partial + \partial^* \bar \partial ( \partial^* \bar \partial )^* + ( \partial^* \bar \partial )^* \partial^* \bar \partial + \partial^* \partial + \bar \partial^* \bar \partial
\]
is a self-adjoint elliptic operator on a compact almost complex manifold endowed with an almost Hermitian metric. Hence, its kernel is a finite dimensional vector space over $\C$, whose dimension $h^{p,q}_{\partial + \bar \partial}$ is independent of the choice of Hermitian metric if $J$ is a complex structure.\\
While in the complex case the main operators that have been studied are $\partial$, $\bar \partial$ (or, equivalently, $d$, $d^c := J^{-1} d J$), in the almost complex case they admit several non-equivalent generalizations, and one should formulate Kodaira-Spencer's question for each one of them (see section \ref{subsec:prel:almostcomplex} for precise definitions of the generalizations which are used in the formulation of the main problem). Therefore, our main problem consists in determining whether or not the dimension of the kernel of Dolbeault-like and Bott-Chern-like Laplacians, built starting from those operators, are independent of the choice of almost Hermitian structure.
\vspace{.2cm}

\textbf{Main problem.} (Generalization of Kodaira-Spencer's problem)\\
Let $(M,J)$ be a compact almost complex $4$-manifold. Fix an almost Hermitian structure and for $P \in \{ \bar \partial, \,  \partial + \bar \partial, \, \bar \delta, \, \delta + \bar \delta, \, d, \, d+d^c\}$ consider the dimension of the spaces of bigraded $P$-harmonic forms $h^{\bullet, \bullet}_P$, if $P \in \{ \partial, \,  \partial + \bar \partial, \, d \}$, and of graded $P$-harmonic forms $h^{\bullet}_P$, if $P \in \{ \bar \delta, \, \delta + \bar \delta, \, d+d^c \}$.
\begin{itemize}
    \item Are $h^{\bullet, \bullet}_P$, $h^{\bullet}_P$ independent of the choice of  almost Hermitian structure?

    \item If $J$ admits a compatible almost K\"ahler structure, are $h^{\bullet, \bullet}_P$, $h^{\bullet}_P$ independent of the choice of almost K\"ahler structure?
\end{itemize}
The same problem can be formulated for manifolds of arbitrary dimension.
\vspace{.2cm}

In this paper we combine new results with results of Cirici and Wilson \cite{CW20a}, Holt \cite{Hol22}, Tardini and the second author \cite{TT20}, and the authors of the present paper \cite{ST23b, ST23c}, in order to prove the following:

\begin{namedtheorem}[Main]
    Let $(M,J)$ be a compact almost complex $4$-manifold. Then:
    \begin{itemize}

        \item the numbers $h^k_{\bar \delta}$, $h^{p,q}_{d}$ depend on the choice of almost Hermitian metric;

        \item the numbers $h^k_{\bar \delta}$, $h^{p,q}_{d}$, $h^k_{\delta + \bar \delta}$, $h^k_{d + d^c}$ do not depend on the choice of almost K\"ahler metric.
    \end{itemize}
\end{namedtheorem}

Our main theorem gives a full answer to the problem we posed for the operators $\delta$, $\bar \delta$ and $d$, $d^c$ on almost K\"ahler $4$-manifolds. It shows that neither $h^k_{\bar \delta}$ nor $h^{p,q}_{d}$ are a metric-independent generalization of Hodge numbers if the metric is not almost K\"ahler. Since none among $h^{p,q}_{\bar \partial}$, $h^k_{\bar \delta}$, $h^{p,q}_d$ is completely independent of the choice of almost Hermitian metric, the original Kodaira-Spencer's problem remains unsolved. The theorem also shows that the numbers $h^k_{\bar \delta}$, $h^{p,q}_{d}$, $h^k_{\delta + \bar \delta}$, $h^k_{d + d^c}$ are independent of the choice of almost K\"ahler metric, and it provides a positive answer to the second question of our main problem. Note that Holt and Zhang \cite{HZ22b} proved that $h^{0,1}_{\bar \partial}$ might depend on the choice of almost K\"ahler metric, while Holt \cite{Hol22} proved that the same is true for $h^{1,2}_{\partial + \bar \partial}$. Here, we prove that for the numbers $h^{p,q}_{\partial + \bar \partial}$ on almost K\"ahler $4$-manifolds, this is the only degree of freedom (Corollary \ref{cor:12}): $h^{p,q}_{\partial + \bar \partial}$ is independent of the choice of almost K\"ahler metric if $(p,q) \neq (1,2)$.
\vspace{.2cm}

As already observed by Cirici and Wilson \cite{CW20b} for $h^{p,q}_d$, it turns out that, on almost K\"ahler $4$-manifolds, several numbers among $h^k_{\bar \delta}$, $h^{p,q}_{d}$, $h^k_{\delta + \bar \delta}$, $h^k_{d + d^c}$, are not only metric-independent, but topological. Hence, it is natural to wonder how many genuine (not topological) almost complex invariants are there among them. As an application of our theory, we find that on almost K\"ahler $4$-manifolds, there are essentially only two almost complex invariants, which do not completely determine each other: the number $h^1_{d+d^c}$ and the number $h^-_J$, the dimension of the space of $J$-anti-invariant harmonic forms (see \cite{DLZ13, DLZ10, LZ09}).

\begin{reptheorem}{thm:invariants}
    Let $(M,J)$ be a compact almost complex $4$-manifold admitting a $J$-compatible almost K\"ahler metric. Then for every choice of $J$-compatible almost K\"ahler metric, the invariants $h^k_{\bar \delta}$, $h^{p,q}_{d}$, $h^k_{\delta + \bar \delta}$, $h^k_{d + d^c}$, are completely determined by:
    \begin{itemize}
        \item the oriented topology of the underlying manifold (more precisely, by the numbers $b_1$ and $b^-$);

        \item the almost complex invariant $h^1_{d+d^c}$;

        \item the almost complex invariant $h^-_J$.
    \end{itemize}
Moreover, the invariants $h^1_{d+d^c}$, $h^-_J$ do not completely determine each other.
\end{reptheorem}

The structure of the paper is as follows: in section \ref{sec:prel}, we briefly recall basic facts on complex and almost complex manifolds, and we review the progresses made on Kodaira-Spencer's problem in the last years. Section \ref{sec:gen:KS} is dedicated to the study of several spaces of harmonic forms. It contains the proof of our main theorem, together with some open questions. Finally, section \ref{sec:applications} contains theorem \ref{thm:invariants} together with other applications of our results.

\section{Preliminaries}\label{sec:prel}

This section begins with a short review of Dolbeault, Bott-Chern and Aeppli cohomologies and the corresponding Laplacians on complex manifolds. Then, we discuss almost complex manifolds and possible generalizations of the mentioned cohomologies. Finally, we focus on Kodaira-Spencer's problem and review recent literature and progresses made on the problem.

\subsection{Complex manifolds.}\label{subsec:prel:complex}

Let $(M,J)$ be a compact complex $2m$-manifold, that is, a compact smooth $2m$-manifold endowed with a smooth $(1,1)$-tensor $J$ such that $J^2 = -\Id$ and satisfying the integrability condition
\[
N_J(X,Y) := [JX,JY] -J[JX,Y]-J[X,JY]-[X,Y] =0
\]
for all $X$,$Y$ vector fields on $M$. Let $A^k = A^k_\C (M)$ be the space of smooth complex $k$-forms on $M$ and let 
\begin{equation}\label{eq:bigraded:dec}
    A^k = \bigoplus_{p+q=k} A^{p,q}
\end{equation}
be the bigraded splitting of $k$-forms induced by $J$. The exterior derivative $d$ decomposes with respect to \eqref{eq:bigraded:dec} as $d = \partial + \bar \partial$, where
\[
\partial := \pi^{p+1,q} \circ d_{|_{A^{p,q}}} \quad \text{and} \quad \bar \partial := \pi^{p,q+1} \circ d_{|_{A^{p,q}}}.
\]
Since $d^2 =0$, the operators $\partial$, $\bar \partial$ satisfy the equations
\begin{equation}
\begin{cases}\label{eq:anticommute}\tag{$\triangle$}
    \partial^2=0, \\
    \bar \partial^2 =0, \\
    \partial \bar \partial + \bar \partial \partial =0,
\end{cases}
\end{equation}
and $(A^{\bullet,\bullet},\partial, \bar \partial)$ is a double complex. The natural cohomologies associated to the double complex are:
\begin{itemize}
    \item the \emph{Dolbeault cohomology}
    \[
    H^{\bullet,\bullet}_{\bar \partial} (M,J) = \frac{\ker \bar \partial \cap A^{\bullet,\bullet}}{\Ima \bar \partial \cap A^{\bullet,\bullet}};
    \]

    \item the \emph{Bott-Chern cohomology} \cite{BC65}
    \[
    H^{\bullet,\bullet}_{BC} (M,J) = \frac{\ker \partial \cap \ker \bar \partial \cap A^{\bullet,\bullet}}{\Ima \partial \bar \partial \cap A^{\bullet,\bullet}};
    \]

    \item the \emph{Aeppli cohomology} \cite{Aep65}
    \[
    H^{\bullet,\bullet}_{A} (M,J) = \frac{\ker \partial \bar \partial \cap A^{\bullet,\bullet}}{(\Ima \partial + \Ima \bar \partial) \cap A^{\bullet,\bullet}}.
    \]
\end{itemize}
Fix a $J$-compatible Hermitian metric on $M$. For any differential operator $P \colon A^\bullet \rightarrow A^\bullet$, denote by $P^*$ its formal adjoint with respect to the metric, and consider the following Laplacians:
\begin{itemize}
    \item the \emph{Dolbeault Laplacian}
    \[
    \Delta_{\bar \partial} = \bar \partial \bar \partial^* + \bar \partial^* \bar \partial;
    \]
    \item the \emph{Bott-Chern Laplacian} \cite{Sch07}
    \[
    \Delta_{BC} := \partial \bar \partial ( \partial \bar \partial)^* + ( \partial \bar \partial)^* \partial \bar \partial + \partial^* \bar \partial ( \partial^* \bar \partial )^* + ( \partial^* \bar \partial )^* \partial^* \bar \partial + \partial^* \partial + \bar \partial^* \bar \partial;
    \]
    \item the \emph{Aeppli Laplacian} \cite{Sch07}
    \[
    \Delta_{A} := \partial \bar \partial ( \partial \bar \partial)^* + ( \partial \bar \partial)^* \partial \bar \partial + \partial^* \bar \partial ( \partial^* \bar \partial )^* + ( \partial^* \bar \partial )^* \partial^* \bar \partial + \partial \partial^* + \bar \partial \bar \partial^*.
    \]
\end{itemize}
The kernel of $\Delta_{\bar \partial}$ (resp.\ $\Delta_{BC}$, $\Delta_A$) restricted to $(p,q)$-forms is denoted by $\H^{p,q}_{\bar \partial}$ (resp.\ $\H^{p,q}_{BC}$, $\H^{p,q}_{A}$) and it is called the \emph{space of Dolbeault-harmonic $(p,q)$-forms} (resp.\ \emph{Bott-Chern-harmonic}, \emph{Aeppli-harmonic}). Since $\Delta_{\bar \partial}$, $\Delta_{BC}$, $\Delta_A$ are self-adjoint elliptic operators on compact manifolds, their kernels are finite dimensional vector spaces over $\C$ and their complex dimensions $h^{p,q}_{\bar \partial}$, $ h^{p,q}_{BC}$, $h^{p,q}_A$, are invariants of the Hermitian structure, which in principle depend on the choice of metric. However, Hodge theory allows to establish that there are isomorphisms
\[
H^{p,q}_{\bar \partial} \cong \H^{p,q}_{\bar \partial}, \quad \quad H^{p,q}_{BC}\cong \H^{p,q}_{BC}, \quad \quad  H^{p,q}_{A} \cong \H^{p,q}_{A},
\]
from which one can conclude that:
\begin{itemize}
    \item Dolbeault, Bott-Chern and Aeppli cohomologies are finite-dimensional vector spaces over $\C$;
    \item the numbers $h^{p,q}_{\bar \partial}$, $ h^{p,q}_{BC}$, $h^{p,q}_A$ are \emph{independent} of the choice of Hermitian structure and are complex invariants.
\end{itemize}

\subsection{Almost complex manifolds.}\label{subsec:prel:almostcomplex}

Let $(M,J)$ be a compact almost complex $2m$-manifold, that is, a compact smooth $2m$-manifold endowed with a smooth $(1,1)$-tensor $J$ such that $J^2 = - \Id$ (with $N_J$ not necessarily vanishing). Existence of $J$ still induces the bigraded decomposition \eqref{eq:bigraded:dec}, but the exterior derivative decomposes as $d = \mu + \partial + \bar \partial + \bar \mu$, with 
\begin{align*}
&\mu := \pi^{p+2,q-1} \circ d_{|_{A^{p,q}}}, \quad \quad   \partial := \pi^{p+1,q} \circ d_{|_{A^{p,q}}}, \\ 
&\bar \partial := \pi^{p,q+1} \circ d_{|_{A^{p,q}}}, \quad \quad \bar \mu := \pi^{p-1,q+2} \circ d_{|_{A^{p,q}}}.
\end{align*}
An almost complex structure is said to be \emph{integrable} if it is a complex structure, i.e., if $N_J=0$. Equivalently, $J$ is integrable if and only if $\bar \mu =0$, if and only if $\bar \partial^2=0$. The equation $d^2=0$ implies that
\begin{equation}
\begin{cases}\label{eq:no:anticommute}\tag{$\triangledown$}
    \bar \mu^2=0, \\
    \bar \mu \bar \partial + \bar \partial \bar \mu  =0, \\
    \bar \partial^2 + \bar \mu \partial + \partial \bar \mu =0, \\
    \mu \bar \mu + \bar \mu \mu + \partial \bar \partial + \bar \partial \partial =0.
\end{cases}
\end{equation}
Since \eqref{eq:anticommute} does not hold, then $(A^{\bullet,\bullet}, \partial, \bar \partial)$ is not a double complex and the usual definitions of Dolbeault, Bott-Chern and Aeppli cohomologies are not well-given for a strict almost complex structure.\\
After observing that by \eqref{eq:no:anticommute} the $\bar \mu$-cohomology is well-defined and that $\bar \partial^2 =0$ on it, Cirici and Wilson gave a definition of \emph{Dolbeault cohomology of an almost complex manifold} obtained as the $\bar \partial$-cohomology computed on the $\bar \mu$-cohomology \cite{CW21}.\\
Later, Coelho, Placini and Stelzig gave a definition of Bott-Chern and Aeppli cohomologies obtained applying the usual definition to suitable subcomplex and quotient complex of $A^{\bullet,\bullet}$ \cite{CPS22}.\\
Inspired by the symplectic cohomologies $H^\bullet_{d +d^\Lambda}$, $H^\bullet_{d d^\Lambda}$, introduced by Tseng and Yau \cite{TY12a, TY12b, TY14}, we gave a definition of Bott-Chern and Aeppli cohomologies $H^\bullet_{d+d^c}$, $H^\bullet_{d d^c}$, based on the operators $d$, $d^c := J^{-1} d J$, rather than $\partial$, $\bar \partial$ \cite{ST23b}. In the complex case the two approaches are equivalent since one has that $d = \partial + \bar \partial$, $d^c = i (\bar \partial - \partial)$, obtaining the equality
\[
H^{\bullet,\bullet}_{BC} = \frac{\ker \partial \cap \ker \bar \partial \cap A^{\bullet,\bullet}}{\Ima \partial \bar \partial \cap A^{\bullet,\bullet}} = \frac{\ker d \cap \ker d^c \cap A^{\bullet,\bullet}}{\Ima d d^c \cap A^{\bullet,\bullet}}
\]
for Bott-Chern cohomology (and a similar equality for Aeppli cohomology). In the almost complex case the approach based on $d$, $d^c$ is not equivalent to the one based on $\partial$, $\bar \partial$. Rather, the operators $\delta := \partial + \bar \mu$, $\bar \delta := \bar \partial + \mu$ appear as an appropriate generalization of $\partial$, $\bar \partial$ to almost complex manifolds, at least in the context of Bott-Chern and Aeppli cohomologies (cf.\ \cite{ST23b}).

\subsection{The Kodaira-Spencer's problem.}\label{subsec:prel:KS}

Fix a Hermitian metric on a compact complex manifold $(M,J)$. Then the Hodge numbers
\[
h^{p,q}_{\bar \partial} := \dim_\C (A^{p,q} \cap \ker \Delta_{\bar \partial})
\]
are independent of the choice of Hermitian structure on $(M,J)$ (see section \ref{subsec:prel:complex}).\\
Suppose now that $(M,J)$ is a compact almost complex manifold and fix an \emph{almost Hermitian metric}, i.e., a $J$-compatible metric inducing a Hermitian metric on the fibers of $TM$. Since the Dolbeault Laplacian 
\[
\Delta_{\bar \partial} = \bar \partial \bar \partial  ^* + \bar \partial  ^* \bar \partial  
\]
is a self-adjoint elliptic operator and the manifold is compact, its kernel $\H^{p,q}_{\bar \partial} = A^{p,q} \cap \ker \Delta_{\bar \partial}$ is a finite-dimensional vector space over $\C$, called the \emph{space of $\bar \partial$-harmonic $(p,q)$-forms}. Denote by $h^{p,q}_{\bar \partial}$ its complex dimension.
\vspace{.2cm}

\textbf{Question.} (Kodaira-Spencer, see \cite[Problem 20]{Hir54})\\
Is $h^{p,q}_{\bar \partial}$ independent of the choice of (almost) Hermitian structure? If not, give some other definition of the $h^{p,q}_{\bar \partial}$ which depends only on the almost complex structure and which generalizes the $h^{p,q}_{\bar \partial}$ of a complex manifold.
\vspace{.2cm}

A negative answer to Kodaira-Spencer's question was given by Holt and Zhang, who proved that there exist almost complex structures on the Kodaira-Thurston manifold such that $h^{0,1}_{\bar \partial}$ varies with different choices of almost Hermitian metric \cite[Theorem 5.1]{HZ22a}.\\
Following the observation made by Hirzebruch on the ellipticity of the Dolbeault Laplacian on almost complex manifolds, it is not hard to verify that most of the Laplacians appearing as a natural generalization of Dolbeault, Bott-Chern and Aeppli Laplacians are still self-adjoint and elliptic, even when $J$ is not integrable. For instance, we can consider a generalization of Bott-Chern Laplacian
\[
\Delta_{BC} := \partial \bar \partial ( \partial \bar \partial)^* + ( \partial \bar \partial)^* \partial \bar \partial + \partial^* \bar \partial ( \partial^* \bar \partial )^* + ( \partial^* \bar \partial )^* \partial^* \bar \partial + \partial^* \partial + \bar \partial^* \bar \partial
\]
which is self-adjoint and elliptic on almost complex manifolds \cite{PT22a}, it coincides with the usual Bott-Chern Laplacian when $J$ is integrable, and the dimension $h^{p,q}_{BC}$ of its kernel is an almost Hermitian invariant. Therefore, the following appears as a natural problem.
\vspace{.2cm}

\textbf{Question.} Is $h^{p,q}_{BC}$ independent of the choice of almost Hermitian structure? If not, give some other definition of the $h^{p,q}_{BC}$ which depends only on the almost complex structure and which generalizes the $h^{p,q}_{BC}$ of a complex manifold.
\vspace{.2cm}

The above question can be considered as a natural generalization of Kodaira-Spencer's question from Hodge numbers to Bott-Chern numbers. Even in this case the answer is negative since $h^{1,2}_{BC}$ varies with different choices of metric \cite{Hol22}. The same considerations and questions are valid for Aeppli cohomology. 

Nevertheless, in the almost complex case there is more than one possible generalization of Bott-Chern Laplacian and approaches based on Laplacians built using different operators are not equivalent (cf.\ \cite[Remark 4.15]{ST23b}). We summarize the directions investigated so far in the literature.
\begin{itemize}
    \item The approach using $\partial$,$\bar \partial$ has been studied by Cattaneo, Tardini and the second author \cite{CTT22}, Holt \cite{Hol22}, Holt and Piovani \cite{HP23}, Holt and Zhang \cite{HZ22b, HZ22a}, Piovani and the second author \cite{PT22b}, Tardini and the second author \cite{TT21, TT22}, for what concerns the Laplacian $\Delta_{\bar \partial}$ and by Holt \cite{Hol22}, Holt and Piovani \cite{HP23}, Piovani and Tardini \cite{PTa23}, Piovani and the second author \cite{PT22a}, for what concerns the Laplacian $\Delta_{BC}$. The related invariants are the numbers $h^{p,q}_{\bar \partial}$ and $h^{p,q}_{BC}$. We denote the latter by $h^{p,q}_{\partial + \bar \partial}$ to avoid confusion with other spaces of Bott-Chern-like harmonic forms.

    \item The approach using $\delta$,$\bar \delta$ has been introduced and studied by Tardini and the second author \cite{TT20}. They considered the Laplacians
    \[
    \Delta_{\bar \delta} := \bar \delta \bar \delta^* + \bar \delta^* \bar \delta 
    \]
    and
    \[
    \Delta_{\delta + \bar \delta} := \delta \bar \delta ( \delta \bar \delta)^* + ( \delta \bar \delta)^* \delta \bar \delta + \delta^* \bar \delta ( \delta^* \bar \delta )^* + ( \delta^* \bar \delta )^* \delta^* \bar \delta + \delta^* \delta + \bar \delta^* \bar \delta.
    \]
    The dimensions of the related spaces of harmonic forms are $h^k_{\bar \delta}$ and $h^k_{\delta + \bar \delta}$.

    \item The approach using $d$,$d^c$ was introduced by Cirici and Wilson for $d$-harmonic $(p,q)$-forms on almost K\"ahler manifolds \cite{CW20b} and was furhter studied by Holt, Piovani and the second author on almost complex manifolds \cite{HPT23}. Later, the authors of the present paper introduced Bott-Chern-like harmonic forms based on $d$,$d^c$, called $(d+d^c)$-harmonic forms \cite{ST23b}. The theory of $(d+d^c)$-harmonic forms is further developed in \cite{ST23c} and in the present paper. In this case, the invariants are $h^{p,q}_d$ and $h^k_{d+d^c}$.
\end{itemize} 

As first observed in \cite{CW20b} (see also \cite{HZ22b, HZ22a}), results on all of the mentioned invariants improve if we assume that $2m=4$ and that the metric is almost K\"ahler. Hence, we formulate a problem which appears as a natural generalization of Kodaira-Spencer's problem.
\vspace{.2cm}

\textbf{Main problem.} (Generalization of Kodaira-Spencer's problem)\\
Let $(M,J)$ be a compact almost complex $4$-manifold. Fix an almost Hermitian structure and for $P \in \{ \bar \partial, \,  \partial + \bar \partial, \, \bar \delta, \, \delta + \bar \delta, \, d, \, d+d^c\}$ consider the dimension of the spaces of harmonic forms $h^{\bullet, \bullet}_P$, if $P \in \{ \partial, \,  \partial + \bar \partial, \, d \}$, or $h^{\bullet}_P$, if $P \in \{ \bar \delta, \, \delta + \bar \delta, \, d+d^c \}$.
\begin{itemize}
    \item Are $h^{\bullet, \bullet}_P$, $h^{\bullet}_P$ independent of the choice of  almost Hermitian structure?

    \item Are $h^{\bullet, \bullet}_P$, $h^{\bullet}_P$ independent of the choice of almost K\"ahler structure?
\end{itemize}
Even though we formulated the problem for $4$-manifolds, it makes sense for manifolds of arbitrary even dimension.
\vspace{.2cm}

For convenience of the reader, we write explicitly the spaces of harmonic forms treated in the present paper. The space of \emph{$\bar \partial$-harmonic $(p,q)$-forms} is
\[
\H^{p,q}_{\bar \partial} = \{ \alpha \in A^{p,q} : \bar \partial \alpha =0, \, \, \partial * \alpha =0 \}.
\]
The space of \emph{$(\partial + \bar \partial)$-harmonic $(p,q)$-forms} is
\[
\H^{p,q}_{\bar \partial} = \{ \alpha \in A^{p,q} : \partial \alpha =0, \, \,  \bar \partial \alpha =0, \, \, \partial \bar \partial * \alpha =0 \}.
\]
The space of \emph{$\bar \delta$-harmonic $k$-forms} is
\[
\H^k_{\bar \delta} = \{ \alpha \in A^k : \bar \delta \alpha =0, \, \, \delta * \alpha =0 \}.
\]
The space of \emph{$(\delta + \bar \delta)$-harmonic $k$-forms} is
\[
\H^k_{\delta + \bar \delta} = \{ \alpha \in A^k :  \delta \alpha =0, \, \, \bar \delta  \alpha =0, \, \, \delta \bar \delta * \alpha =0 \}.
\]
The space of \emph{$d$-harmonic $(p,q)$-forms} is
\[
\H^{p,q}_{d} = \{ \alpha \in A^{p,q} : d \alpha =0, \, \, d * \alpha =0 \}.
\]
The space of \emph{$(d + d^c)$-harmonic $k$-forms} is
\[
\H^k_{d + d^c} = \{ \alpha \in A^k :  d \alpha =0, \, \, d^c  \alpha =0, \, \, d^c d * \alpha =0 \}.
\]

\section{Generalizations of the Kodaira-Spencer's problem}\label{sec:gen:KS}

In this section we study the metric-independence of the numbers $h^k_{\bar \delta}$, $h^k_{\delta + \bar \delta}$, $h^{p,q}_d$, $h^{p,q}_{d+d^c}$ on compact almost Hermitian $4$-manifolds and we prove our main theorem.
\vspace{.2cm}

We begin by proving a decomposition for $\bar \delta$ and $(\delta + \bar \delta)$-harmonic $2$-forms. Denote by $\H^{(2,0)(0,2)}_J$ the spaces of $J$-anti-invariant, $d$-harmonic $2$-forms
\[
\H^{(2,0)(0,2)}_J = \{ \alpha \in A^{2,0} \oplus A^{0,2} : d \alpha =0, \, \, d*\alpha =0 \}.
\]

\begin{theorem}\label{thm:decomposition}
    Let $(M,J)$ be a compact almost complex $4$-manifold. Fix a $J$-compatible almost Hermitian metric on $M$. Then:
    \begin{itemize}
        \item there is a decomposition of $\bar \delta$-harmonic $2$-forms
        \[
        \H^2_{\bar \delta} = \H^{1,1}_{\bar \partial} \oplus \H^{(2,0)(0,2)}_J;
        \]
        
        \item there is a decomposition of $(\delta + \bar \delta)$-harmonic $2$-forms
        \[
        \H^2_{\delta + \bar \delta} = \H^{1,1}_{\partial + \bar \partial} \oplus \H^{(2,0)(0,2)}_J.
        \]
    \end{itemize}
\end{theorem}
\begin{proof}
    It is immediate to verify that the following inclusions hold:
    \[
    \H^{1,1}_{\bar \partial} \oplus \H^{(2,0)(0,2)}_J \subseteq \H^2_{\bar \delta}, \quad \quad \quad \H^{1,1}_{\partial + \bar \partial} \oplus \H^{(2,0)(0,2)}_J \subseteq \H^2_{\delta + \bar \delta}.
    \]
    We first prove that $\H^2_{\bar \delta} \subseteq \H^{1,1}_{\bar \partial} \oplus \H^{(2,0)(0,2)}_J$. Let $\alpha \in \H^2_{\bar \delta}$. Then $\bar \delta \alpha =0$ and $\delta * \alpha =0$. Writing $\alpha$ as the sum of bigraded forms, $\alpha = \alpha^{2,0} + \alpha^{1,1} + \alpha^{0,2}$, and imposing that $\alpha$ is $\bar \delta$-harmonic, we have that
    \begin{align*}
    0 &= \bar \delta \alpha = (\bar \partial + \mu )( \alpha^{2,0} + \alpha^{1,1} + \alpha^{0,2} ) = \\
    &= \bar \partial \alpha^{2,0} + \bar \partial \alpha^{1,1} + \bar \partial \alpha^{0,2} + \mu \alpha^{2,0} + \mu \alpha^{1,1} + \mu \alpha^{0,2}=\\
    &= \bar \partial \alpha^{2,0} + \mu \alpha^{0,2} + \bar \partial \alpha^{1,1},
    \end{align*}
    since several terms vanish by bidegree reasons. Separating the bidegrees, we obtain that
    \begin{equation}\label{eq:bardelta}
    \bar \partial \alpha^{1,1} = 0 \quad \text{and} \quad  \bar \partial \alpha^{2,0} + \mu \alpha^{0,2} =0.
    \end{equation}
    Similarly, since $(2,0)$ and $(0,2)$-forms are self-dual, we see that 
    \begin{align*}
    0 &= \delta * \alpha = (\partial + \bar \mu ) * ( \alpha^{2,0} + \alpha^{1,1} + \alpha^{0,2} ) = \\
    &= (\partial + \bar \mu )( \alpha^{2,0} + * \alpha^{1,1} + \alpha^{0,2} ) = \\
    &= \partial \alpha^{0,2} + \bar \mu \alpha^{2,0} + \partial * \alpha^{1,1},
    \end{align*}
    from which we get the equations 
    \begin{equation}\label{eq:delta*}
    \partial * \alpha^{1,1} = 0 \quad \text{and} \quad \partial \alpha^{0,2} + \bar \mu \alpha^{2,0} =0.
    \end{equation}
    Combining \eqref{eq:bardelta} and \eqref{eq:delta*}, we immediately deduce that $\alpha^{1,1} \in \H^{1,1}_{\bar \partial}$ and that $\alpha^{2,0} + \alpha^{0,2}$ is $d$-closed and self-dual, thus $d$-harmonic, proving the first part of the theorem. For the second part, we prove the inclusion $\H^2_{\delta + \bar \delta} \subseteq \H^{1,1}_{\partial + \bar \partial} \oplus \H^{(2,0)(0,2)}_J$. Let $ \alpha \in \H^2_{\delta + \bar \delta}$. Then $\delta \alpha =0$, $\bar \delta \alpha =0$ and $\delta \bar \delta * \alpha =0$. Writing $\alpha$ as the sum of forms of pure bidegree, following the same computations of the first part of the proof of this theorem and imposing the conditions $\delta \alpha =0$, $\bar \delta \alpha =0$, we get that
    \begin{equation}\label{eq:delta+bardelta}
    \begin{split}
    &\partial \alpha^{1,1} = 0, \quad (\partial + \bar \mu ) (\alpha^{2,0}+\alpha^{0,2})=0, \\
    &\bar \partial \alpha^{1,1} =0, \quad (\bar \partial + \mu ) (\alpha^{2,0}+\alpha^{0,2})=0.
    \end{split}
    \end{equation}
    For the equation $\delta \bar \delta * \alpha =0$, observe that in general
    \begin{equation}
    \delta \bar \delta = (\partial + \bar \mu)(\bar \partial + \mu) = \partial \bar \partial + \partial \mu + \bar \mu \bar \partial + \bar \mu \mu,
    \end{equation}
    while on $4$-manifolds the equation simplifies to $\delta \bar \delta = \partial \bar \partial + \bar \mu \mu$ and the operator $\delta \bar \delta$ has bidegree $(1,1)$. Thus, by bidegree, we have that
    \begin{equation}\label{eq:deltabardelta}
        0 = \delta \bar \delta * \alpha = (\partial \bar \partial + \bar \mu \mu) * \alpha^{1,1} = \partial \bar \partial * \alpha^{1,1}.
    \end{equation}
    Finally, from \eqref{eq:delta+bardelta} and \eqref{eq:deltabardelta}, we conclude that $\alpha^{1,1} \in \H^{1,1}_{\partial + \bar \partial}$ and $\alpha^{2,0}+\alpha^{0,2} \in \H^{(2,0)(0,2)}_J$. This completes the proof of the theorem.
\end{proof}

\begin{cor}\label{cor:bardelta}
    Let $(M,J)$ be a compact almost complex $4$-manifold. Fix a $J$-compatible almost Hermitian metric on $M$. Then $h^2_{\delta + \bar \delta} = b^-+1 + h^-_J$ and it is independent of the choice of almost Hermitian metric. 
    If the metric is almost K\"ahler, then $h^2_{\bar \delta} = h^2_{\delta + \bar \delta}$ and it is independent of the choice of almost K\"ahler metric.
\end{cor}
\begin{proof}
    By theorem \ref{thm:decomposition}, we have that $h^2_{\delta + \bar \delta} = h^{1,1}_{\partial + \bar \partial} + h^-_J$. By \cite[Theorem 4.2]{Hol22}, on every almost Hermitian $4$-manifold we have $h^{1,1}_{\partial + \bar \partial} = b^-+1$. This proves the first part of the corollary. The second part follows either using the equality $h^k_{\bar \delta} = h^k_{\delta + \bar \delta}$, valid on almost K\"ahler manifolds \cite[Proposition 6.10]{TT20}, or observing that by theorem \ref{thm:decomposition} we have $h^2_{\bar \delta} = h^{1,1}_{\bar \partial}  + h^-_J$ and by \cite[Proposition 6.1]{HZ22a}, $h^{1,1}_{\bar \partial} = b^- +1$ if the metric is almost K\"ahler.
\end{proof}

We are ready to prove our main theorem.

\begin{proof}[Proof of main theorem]
First, we treat the almost K\"ahler case. By \cite[Proposition 6.10]{TT20}, on almost K\"ahler manifolds we have that $\H^k_{\bar \delta} = \H^k_{\delta + \bar \delta}$. Therefore, we just need to prove the theorem for $h^k_{\bar \delta} $ or $h^k_{\delta + \bar \delta}$. If $k = 0,4$, one easily checks that the only $\bar \delta$-harmonic functions are constant and that the only $\bar \delta$-harmonic $4$-forms are constant multiple of the volume form, therefore $h^0_{\bar \delta} = h^4_{\bar \delta} = 1$. By corollary \ref{cor:bardelta}, $h^2_{\bar \delta}$ is independent of the choice of almost K\"ahler metric.  Finally, we have that $h^1_{\delta + \bar \delta}$ is independent of the choice of almost Hermitian metric since $\alpha \in \H^1_{\delta + \bar \delta}$ if and only if $\delta \alpha = \bar \delta \alpha =0$. Using the equality $h^k_{\bar \delta} = h^k_{\delta + \bar \delta}$, valid for almost K\"ahler metrics, and the equality $h^k_{\bar \delta} = h^{4-k}_{\bar \delta}$, valid for almost Hermitian metrics \cite[Remark 5.7]{TT20}, we deduce that
    \[
    h^1_{\delta + \bar \delta} = h^1_{\bar \delta} = h^3_{\bar \delta} = h^3_{\delta + \bar \delta},
    \]
    proving metric-independence of $h^k_{\bar \delta} = h^k_{\delta + \bar \delta}$. The metric-independence of the numbers $h^{p,q}_d$ is proved in \cite[Corollary 5.9]{CW20b}. For the metric-independence of $h^k_{d+d^c}$, the cases $k=0, 1,4$ are, again, immediate. The cases $k=2,3$ follow from \cite[Corollary 5.3]{ST23b} and \cite[Corollary 5.5]{ST23b} respectively, after observing that on almost K\"ahler $4$-manifolds $h^3_{d+d^c} = h^3_{d+d^\Lambda} = b_1$ is metric-independent.\\
    To conclude the proof of the theorem, we show that $h^k_{\bar \delta}$ and $h^{p,q}_d$ depend on the choice of almost Hermitian metric. Note that by theorem \ref{thm:decomposition}, $h^2_{\bar \delta} = h^{1,1}_{\bar \partial} + h^-_J$, which is metric-dependent since $h^{1,1}_{\bar \partial}$ depends on the choice of metric \cite[Theorem 3.7]{TT21}. Finally, also $h^{1,0}_d$ depends on the choice of metric by \cite[Theorem 4.1]{HPT23}.
\end{proof}

\begin{remark}
    Note that as a consequence of the proof of the main theorem, we have that, on almost complex $4$-manifolds, $h^k_{\delta + \bar \delta}$ is metric-independent for $k=0,1,2,4$. In \cite{ST23c}, we prove that also $h^k_{d+d^c}$ is metric-independent for $k = 0,1,2,4$. Hence, we ask the following

    \begin{question} 
    Let $(M,J)$ be a compact almost complex $4$-manifolds. Are $h^3_{\delta + \bar \delta}$, $h^3_{d +d^c}$ independent of the choice of almost Hermitian metric?
    \end{question}

    An affirmative answer would provide a metric-independent generalization of Bott-Chern numbers to almost complex $4$-manifolds.
\end{remark}

\section{Applications}\label{sec:applications}
As an application of the theory of section \ref{sec:gen:KS}, in this section we show the following

\begin{theorem}\label{thm:invariants}
    Let $(M,J)$ be a compact almost complex $4$-manifold admitting a $J$-compatible almost K\"ahler metric. Then for every choice of $J$-compatible almost K\"ahler metric, the invariants $h^k_{\bar \delta}$, $h^{p,q}_{d}$, $h^k_{\delta + \bar \delta}$, $h^k_{d + d^c}$, are completely determined by:
    \begin{itemize}
        \item the oriented topology of the underlying manifold (more precisely, by the numbers $b_1$ and $b^-$);

        \item the almost complex invariant $h^1_{d+d^c}$;

        \item the almost complex invariant $h^-_J$.
    \end{itemize}
Moreover, the invariants $h^1_{d+d^c}$, $h^-_J$ do not completely determine each other.
\end{theorem}

In the proof of the theorem, we also explicitly compute $h^k_{\bar \delta}$, $h^{p,q}_{d}$, $h^k_{\delta + \bar \delta}$, $h^k_{d + d^c}$ in terms of $b_1$, $b^-$, $h^1_{d+d^c}$ and $h^-_J$ (see table \ref{tab:main}).

\begin{proof}[Proof of theorem \ref{thm:invariants}]
Fix an arbitrary $J$-compatible almost K\"ahler metric. The cases $k=0,4$ are easy to deal with. For the remaining values of $k$, first note that by \cite[Proposition 6.10]{TT22}, we have that $h^k_{\bar \delta} = h^k_{\delta + \bar \delta}$. By theorem \ref{thm:decomposition} $h^2_{\delta + \bar \delta} = b^- + 1 + h^-_J$ and $h^1_{\delta + \bar \delta} = h^3_{\delta + \bar \delta}$. Hence, the only degrees of freedom for $h^k_{\bar \delta}$, $h^k_{\delta + \bar \delta}$ are $h^1_{\delta + \bar \delta}$ and $h^-_J$.\\
By \cite[Corollary 5.2]{ST23b}, we have that $h^2_{d+d^c} = b^- +1 + h^-_J$ and by \cite[Corollary 5.5]{ST23b} we obtain $h^3_{d+d^c} = b_1$. Again, the only degree of freedom for $h^k_{d+d^c}$ are $h^1_{d + d^c}$ and $h^-_J$. Furthermore, we observe that since $d = \delta + \bar \delta$ and $d^c = i (\bar \delta - \delta)$, one immediately deduces that $h^1_{\delta + \bar \delta} = h^1_{d + d^c}$.\\
For the numbers $h^{p,q}_d$, it is proved in \cite{CW20b} that $h^{2,0}_d = h^{0,2}_d = 0$, that $h^{1,1}_d = b^-+1$, and that $h^{1,0}_d = h^{0,1}_d = h^{2,1}_d = h^{1,2}_d = h^{1,0}_{\bar \partial}$. To conclude the proof of the first part of the theorem, we prove that on compact almost K\"ahler $4$-manifolds $h^1_{d+d^c} = 2 h^{1,0}_d$. By definition (see also \cite{ST23b}), one has that
\[
\H^1_{d + d^c} = A^1 \cap \ker d \cap \ker d^c = (A^{1,0} \cap \ker d) \oplus (A^{0,1} \cap \ker d),
\]
hence $h^1_{d+d^c} = 2 \dim_\C ( A^{1,0} \cap \ker d)$. Again by definition, we also have that $h^{1,0}_d = \dim_\C (A^{1,0} \cap \ker d \cap \ker d^*)$. Finally, note that
\[
A^{1,0} \cap \ker d \cap \ker d^* \subseteq A^{1,0} \cap \ker d \subseteq A^{1,0} \cap \ker \bar \partial = A^{1,0} \cap \ker d \cap \ker d^*,
\]
where the last equality follows from \cite[Proof of corollary 5.9]{CW20b}. This implies the equality of spaces
\[
\H^{1,0}_d = \H^{1,0}_{\bar \partial} = A^{1,0} \cap \ker d
\]
and shows that $h^1_{d+d^c} = 2 h^{1,0}_{d}$. The second part of the theorem follows from the fact that there exists a symplectic $4$-manifold $(M,\omega)$ and a curve of almost complex structures $J_t$, with $t \in (-\epsilon, \epsilon)$, such that:
\begin{itemize}
    \item $\omega$ is an almost K\"ahler metric for each $J_t$;

    \item $h^1_{d+d^c} (J_t)$ varies for different values of $t$;

    \item $h^-_{J_t} =0$ for all $t \in (-\epsilon, \epsilon)$.
\end{itemize}
The symplectic $4$-manifold and the curve of almost complex structures that we have to consider are those given in \cite[Example 6.5]{ST23b}, where the first two claims we made on $J_t$ are also proved. To prove that $h^-_{J_t} =0$ for all $t \in (-\epsilon, \epsilon)$, note that for the considered manifold we have $b^+=b^-=1$ and that by \cite[Corollary 3.4]{DLZ10}, if $b^+=1$ then $h^-_J =0$ for all tamed almost complex structures.
\end{proof}

In the proof of theorem \ref{thm:invariants}, we used the fact, proved by Cirici and Wilson \cite{CW20b}, that on almost K\"ahler $4$-manifolds we have 
\[
A^{1,0} \cap \ker d \cap \ker d^* = A^{1,0} \cap \ker \bar \partial.
\]
Their proof involves the use of the almost K\"ahler identities. We want to give here a direct proof, based on integration, of the inclusion 
\[
A^{1,0} \cap \ker \bar \partial \subseteq A^{1,0} \cap \ker d \cap \ker d^*.
\]
The opposite inclusion is immediate.

\begin{lemma}\label{lemma:integration}
    Let $(M,J)$ be a compact almost complex $4$-manifold endowed with an almost K\"ahler metric. Then 
    \[
    A^{1,0} \cap \ker \bar \partial \subseteq A^{1,0} \cap \ker d \cap \ker d^*.
    \]
\end{lemma}
\begin{proof}
    Note that a $(1,0)$-form $\alpha$ is $d$-harmonic if and only if $\bar \partial \alpha = \partial \alpha = \bar \mu \alpha = \partial^* \alpha =0$. Indeed, the remaining equations are automatically satisfied by bidegree reasons. Suppose that $\bar \partial \alpha =0$. We have to prove that $ \partial \alpha = \bar \mu \alpha = \partial^* \alpha =0$. Observe that $\partial \alpha$ has bidegree $(2,0)$, while $\bar \mu \alpha$ has bidegree $(0,2)$. Hence $\delta \alpha =0$ if and only if $\partial \alpha =0$ and $\bar \mu \alpha =0$. We have that
    \[
    \norm{\delta \alpha }^2 = \int_M \delta \alpha \wedge \overline{* \delta \alpha} = \int_M (\delta + \bar \delta) \alpha \wedge * (\delta + \bar \delta) \bar \alpha = \int_M d\alpha \wedge * d\bar \alpha,
    \]
    where in the second equality we used the fact that $\bar \delta \alpha = \bar \partial \alpha =0$. Since $d \bar \alpha = \bar \partial \bar \alpha + \mu \bar \alpha$ has bidegree $(2,0) + (0,2)$ for any choice of almost Hermitian metric, $d \bar \alpha$ is necessarily a self-dual form, and we have
    \[
    \norm{\delta \alpha }^2 = \int_M d\alpha \wedge * d\bar \alpha = \int_M d\alpha \wedge d\bar \alpha = \int_M d (\alpha \wedge d\bar \alpha) = 0,
    \]
    by Stokes' theorem, showing that $\partial \alpha = \bar \mu \alpha =0$. For the last equation, we have that
    \[
    \partial^* \alpha = - * \bar \partial * \alpha = i * \bar \partial ( \omega \wedge \alpha),
    \]
    where $\omega$ is the fundamental form of the almost K\"ahler metric. Finally, we compute that
    \[
    \bar \partial (\omega \wedge \alpha) = \bar \partial \omega \wedge \alpha + \omega \wedge \bar \partial \alpha =0,
    \]
    completing the proof of the lemma. Note that this last equation is the only instance in the proof where we use that the metric is almost K\"ahler.
\end{proof}

\begin{remark}
    The first part of the proof of lemma \ref{lemma:integration} also appears in \cite[Lemma 4.1]{Lin23}.
\end{remark}

We conclude the section studying the numbers $h^{p,q}_{\bar \partial}$, $h^{p,q}_{\partial + \bar \partial}$ of compact almost K\"ahler $4$-manifolds.

\begin{lemma}\label{lemma:delbar}
    Let $(M,J)$ be a compact almost complex $4$-manifold admitting a $J$-compatible almost K\"ahler metric. Then for every choice of $J$-compatible almost K\"ahler metric we have that
    \begin{itemize}
        \item [(i)] $\H^{2,1}_{\partial + \bar \partial} \cong \H^{1,0}_{ \bar \partial}$;

        \item [(ii)] $\H^{1,2}_{\partial + \bar \partial} \cong \H^{0,1}_{ \bar \partial}$.
    \end{itemize}
\end{lemma}
\begin{proof}
    The proofs of (i) and (ii) are similar, hence we prove only (i). The isomorphism between $\H^{2,1}_{\partial + \bar \partial}$ and $\H^{1,0}_{ \bar \partial}$ is provided by the $\C$-linear Hodge $*$ operator. Let $\alpha^{1,0} \in \H^{1,0}_{ \bar \partial}$. Then $\bar \partial \alpha =0$. After taking the Hodge $*$, we obtain that $\partial * \alpha^{1,0} = 0$ by bidegree reasons, that $\bar \partial * \alpha^{1,0}$ is proportional to $\bar \partial (\omega \wedge \alpha^{1,0}) = \omega \wedge \bar \partial \alpha^{1,0} = 0$ and that $\partial \bar \partial * (* \alpha^{1,0})=0$. This proves the inclusion 
    \[
    *(\H^{1,0}_{ \bar \partial}) \subseteq \H^{2,1}_{\partial + \bar \partial}.
    \]
    For the opposite inclusion, we need to use the almost K\"ahler identities \cite{CW20b}. Since the Hodge $*$ is an isomorphism between $(2,1)$-forms and $(1,0)$-forms, any form in $\H^{2,1}_{\partial + \bar \partial}$ can be written as $\omega \wedge \alpha^{1,0}$ for some $\alpha^{1,0} \in A^{1,0}$. Moreover, since $\bar \partial (\omega \wedge \alpha^{1,0}) = 0$, the form $\bar \partial \alpha^{1,0}$ is primitive, and since $\partial \bar \partial * (\omega \wedge \alpha^{1,0}) =0$, we have $\partial \bar \partial \alpha^{1,0}=0$. By the almost K\"ahler identities, we have that
    \[
    i \bar \partial^* \bar \partial \alpha^{1,0} = [\Lambda, \partial] \bar \partial \alpha^{1,0} = \Lambda \partial \bar \partial \alpha^{1,0} - \partial \Lambda \bar \partial \alpha^{1,0} =0,
    \]
    where in the last equality we used the fact that $\bar \partial \alpha^{1,0}$ is primitive and $\partial \bar \partial \alpha^{1,0}=0$. This completes the proof. 
\end{proof}

\begin{lemma}\label{lemma:h3}
    Let $(M,J)$ be a compact almost complex $4$-manifold. Then for every choice of $J$-compatible Hermitian metric we have that
    \[
    \H^{3}_{\delta + \bar \delta} = \H^{2,1}_{\partial + \bar \partial} \oplus \overline{\H^{2,1}_{\partial + \bar \partial}}.
    \]
\end{lemma}
\begin{proof}
    Let $\alpha \in \H^{2,1}_{\partial + \bar \partial}$. Then $\bar \partial \alpha =0$ and $\partial \bar \partial * \alpha =0$. Observe that 
\[
\delta \alpha = (\partial + \bar \mu) \alpha =0
\]
by bidegree reasons and that
\[
\bar \delta \alpha = (\bar \partial + \mu ) \alpha = \bar \partial \alpha =0
\]
by bidegree reasons and $\bar \partial \alpha =0$. By \eqref{eq:deltabardelta}, we have that
\[
\delta \bar \delta * \alpha = (\partial \bar \partial + \bar \mu \mu) * \alpha = \partial \bar \partial * \alpha =0
\]
since $\mu \alpha =0$. This shows the inclusion $\H^{2,1}_{\partial + \bar \partial} \subseteq \H^3_{\delta + \bar \delta}$. Noting that the equations $\delta \alpha =0$, $\bar \delta \alpha =0$, $\delta \bar \delta * \alpha =0$ are symmetric by complex conjugation, we also have $\overline{\H^{2,1}_{\partial + \bar \partial}} \subseteq \H^3_{\delta + \bar \delta}$. For the opposite inclusion
\[
\H^3_{\delta + \bar \delta} \subseteq \H^{2,1}_{\partial + \bar \partial} \oplus \overline{\H^{2,1}_{\partial + \bar \partial}},
\]
let $\alpha \in \H^3_{\delta + \bar \delta}$. Write $\alpha$ as the sum of bigraded forms $\alpha = \alpha^{2,1} + \alpha^{1,2}$. By bidegree reasons and the equation $\delta \alpha =0$, we have
\[
0 = \delta \alpha = (\partial + \bar \mu) (\alpha^{2,1} + \alpha^{1,2}) = \partial \alpha^{1,2}.
\]
Similarly, from the equation $\bar \delta \alpha =0$, we deduce that $\bar \partial \alpha^{2,1} = 0$. Finally, from the equation $\delta \bar \delta * \alpha =0$, we get
\[
0 = \delta \bar \delta * \alpha = (\partial \bar \partial + \bar \mu \mu) * (\alpha^{2,1} + \alpha^{1,2}).
\]
Since $\delta \bar \delta$ has bidegree $(1,1)$, we can separate the bidegrees to get two equations
\[
\begin{cases}
    \partial \bar \partial + \bar \mu \mu * \alpha^{2,1} =0,\\
    \partial \bar \partial + \bar \mu \mu * \alpha^{1,2} = 0.
\end{cases}
\]
Observing that $\bar \mu \mu * \alpha^{2,1} =0$ (for bidegree reasons), that $\mu \bar \mu * \alpha^{1,2}=0$ (bidegree reasons) and that $\partial \bar \partial + \bar \mu \mu = - \bar \partial \partial - \mu \bar \mu$, all of our equations reduce to
\[
\begin{cases}
    \bar \partial \alpha^{2,1} = 0,\\
    \partial \bar \partial * \alpha^{2,1} = 0,
\end{cases}
\quad \quad 
\begin{cases}
    \partial \alpha^{1,2} = 0,\\
    \bar \partial \partial * \alpha^{1,2} = 0,
\end{cases}
\]
proving that $\alpha^{2,1} \in \H^{2,1}_{\partial + \bar \partial}$ and $\overline{\alpha^{1,2}} \in \H^{2,1}_{\partial + \bar \partial}$, and thus our lemma.
\end{proof}

\begin{cor}\label{cor:12}
    Let $(M,J)$ be a compact almost complex manifold admitting a $J$-compatible almost K\"ahler metric. Then $h^{p,q}_{\partial + \bar \partial}$ is independent of the choice of $J$-compatible almost K\"alher metric if $(p,q) \neq (1,2)$.
\end{cor}

\begin{proof}
    The result is mostly contained in \cite{Hol22}, except for bidegree $(2,1)$. By lemma \ref{lemma:h3}, we have that $h^{2,1}_{\partial + \bar \partial} = \frac{1}{2} h^3_{\delta + \bar \delta}$, which is metric independent on almost K\"ahler $4$-manifolds by the main theorem.
\end{proof}

\begin{table}[H]
\begin{center}
    \caption{\label{tab:main}The numbers $h^k_{\bar \delta}$, $h^{p,q}_{d}$, $h^k_{\delta + \bar \delta}$, $h^k_{d + d^c}$ of compact almost K\"ahler $4$-manifolds.}
    \end{center}
    
    \begin{adjustbox}{center}
    \begin{tabular}{ c|c c|c c c|c c}
     $k$ & \multicolumn{2}{c|}{$1$} & \multicolumn{3}{c|}{$2$} & \multicolumn{2}{c}{$3$} \\
    $(p,q)$ & $(1,0)$ & $(0,1)$ & $(2,0)$ & $(1,1)$ & $(0,2)$ & $(2,1)$ & $(1,2)$ \\
    \hline
         & & & & & & & \\[-2.3ex]
    $h^k_{\bar \delta}$ & \multicolumn{2}{c|}{$h^1_{d+d^c}$} & \multicolumn{3}{c|}{$b^- + 1 + h^-_J$} & \multicolumn{2}{c}{$h^1_{d+d^c}$} \\
         & & & & & & & \\[-2ex]
    $h^k_{\delta + \bar \delta}$ & \multicolumn{2}{c|}{$h^1_{d+d^c}$} & \multicolumn{3}{c|}{$b^- + 1 + h^-_J$} & \multicolumn{2}{c}{$h^1_{d+d^c}$} \\
         & & & & & & & \\[-2ex]
    $h^{p,q}_{d}$ & $\frac{1}{2}h^1_{d+d^c}$ & $\frac{1}{2}h^1_{d+d^c}$ & $0$ & $b^-+1$ & $0$ & $\frac{1}{2}h^1_{d+d^c}$ & $\frac{1}{2}h^1_{d+d^c}$ \\
         & & & & & & & \\[-2ex]
    $h^k_{d + d^c}$ & \multicolumn{2}{c|}{$h^1_{d+d^c}$} & \multicolumn{3}{c|}{$b^- +1 + h^-_J$} & \multicolumn{2}{c}{$b_1$} \\
    \end{tabular}
    \end{adjustbox}
    \end{table}

{\small
\printbibliography
}

\end{document}